\documentclass[twoside]{article} 
\usepackage{graphicx}
\date{} 
\title{The asymptotic expansion of a sum appearing in an approximate functional 
equation for the Riemann zeta function}
\author{\sc R. B.\ Paris \\
{\em Division of Computing and Mathematics,} \\
{\em Abertay University, Dundee DD1 1HG, UK}}
\begin{document}
\def\f#1#2{\mbox{${\textstyle \frac{#1}{#2}}$}}
\def\dfrac#1#2{\displaystyle{\frac{#1}{#2}}}
\def\boldal{\mbox{\boldmath $\alpha$}}
\newcommand{\bee}{\begin{equation}}
\newcommand{\ee}{\end{equation}}
\newcommand{\sa}{\sigma}
\newcommand{\ka}{\kappa}
\newcommand{\al}{\alpha}
\newcommand{\la}{\lambda}
\newcommand{\ga}{\gamma}
\newcommand{\eps}{\epsilon}
\newcommand{\om}{\omega}
\newcommand{\fr}{\frac{1}{2}}
\newcommand{\fs}{\f{1}{2}}
\newcommand{\g}{\Gamma}
\newcommand{\br}{\biggr}
\newcommand{\bl}{\biggl}
\newcommand{\ra}{\rightarrow}
\newcommand{\gtwid}{\raisebox{-.8ex}{\mbox{$\stackrel{\textstyle >}{\sim}$}}}
\newcommand{\ltwid}{\raisebox{-.8ex}{\mbox{$\stackrel{\textstyle <}{\sim}$}}}
\renewcommand{\topfraction}{0.9}
\renewcommand{\bottomfraction}{0.9}
\renewcommand{\textfraction}{0.05}
\newcommand{\mcol}{\multicolumn}
\date{}
\maketitle
\pagestyle{myheadings}
\markboth{\hfill \sc R. B.\ Paris  \hfill}
{\hfill \sc Asymptotics of a finite sum\hfill}
\begin{abstract}
A representation for the Riemann zeta function valid for arbitrary complex $s=\sa+it$ is $\zeta(s)=\sum_{n=0}^\infty A(n,s)$, where
 \[A(n,s)=\frac{2^{-n-1}}{1-2^{1-s}} \sum_{k=0}^n \bl(\!\!\!\begin{array}{c}n\\k\end{array}\!\!\!\br) \frac{(-)^k}{(k+1)^s}.\]
In this note we examine the asymptotics of $A(n,s)$ as $n\to\infty$ when  $t=an$, where $a>0$ is a fixed parameter, by application of the method of steepest descents to an integral representation.
Numerical results are presented to illustrate the accuracy of the expansion obtained.
\vspace{0.3cm}

\noindent {\bf Mathematics subject classification (2020):} 11M06, 30E15, 34E05, 41A60 
\vspace{0.1cm}
 
\noindent {\bf Keywords:} Riemann zeta function, asymptotic expansions, method of steepest descents, Stokes phenomenon
\end{abstract}

\vspace{0.3cm}

\noindent $\,$\hrulefill $\,$

\vspace{0.3cm}
\begin{center}
{\bf 1.\ Introduction}
\end{center}
\setcounter{section}{1}
\setcounter{equation}{0}
\renewcommand{\theequation}{\arabic{section}.\arabic{equation}}
In 1994, Sondow \cite{Sond} showed by use of Euler's transfomation of series that the Riemann zeta function can be expressed as
\bee\label{e11}
\zeta(s)=\sum_{n=0}^\infty A(n,s)
\ee
for any $s=\sigma+it\in {\bf C}$, where
\bee\label{e12}
A(n,s):=\frac{2^{-n-1}}{1-2^{1-s}} \sum_{k=0}^n \bl(\!\!\!\begin{array}{c}n\\k\end{array}\!\!\!\br) \frac{(-)^k}{(k+1)^s}.
\ee
This formula also appears in an appendix in a paper by Hasse \cite{Hasse} published in 1930.

In a recent paper, Jerby \cite{Jerby} has shown using the Euler-Maclaurin summation formula combined with saddle-point estimates that an approximate functional equation for $\zeta(s)$ is given by
\[\zeta(s)=\sum_{n=0}^{[|t|/\pi N]} A(n,s)+\frac{\chi(s)}{1-2^{s-1}} \sum_{k=1}^N (2k-1)^{s-1}+O(e^{-\om(N)|t|})\]
for positive integer $N$, where
\[\chi(s):=2^s\pi^{s-1} \sin \fs\pi s \,\g(1-s)\]
and $\om(N)>0$. The essential feature of this result is that the error term is exponentially small for large $|t|$. This is in contrast to the classical approximate functional equation given by \cite[p.~97]{Ivic}
\[\zeta(s) =\sum_{n\leq x}n^{-s}+\chi(s) \sum_{n\leq y} n^{s-1}+O(x^{-\sigma})+O(t^{\frac{1}{2}-\sigma} y^{\sigma-1})\]
for $0\leq\sigma\leq 1$ and $x$, $y\geq 1$ with $|t|=2\pi xy$, where the error term is of algebraic order.

An analysis of $A(n,s)$ for large $t>0$ and $n=t/(\pi N)$, $N\geq 1$, was discussed by Jerby in \cite[Section 2]{Jerby} by expressing the sum as an integral over the interval $[0,\infty)$. The saddle-point method was then applied to estimate this integral to show that it is exponentially small. However, no attempt at discussion of the associated paths of steepest descent was made
to establish how the saddles considered were connected to the integration path $[0,\infty)$.

In this note our aim is to examine in some detail the asymptotic expansion of $A(n,s)$ for $n\to\infty$ with $t=an$, where $a>0$ is a fixed parameter. We also employ the method of steepest descents to the same integral representation
of $A(n,s)$ used by Jerby, but carefully consider the connection of the steepest descent paths through certain saddle points with the integration path $[0,\infty)$. 
\vspace{0.6cm}

\begin{center}
{\bf 2.\ The asymptotic expansion of $A(n,s)$}
\end{center}
\setcounter{section}{2}
\setcounter{equation}{0}
\renewcommand{\theequation}{\arabic{section}.\arabic{equation}}
We first express $A(n,s)$ as an integral by noting from (\ref{e12}) that
\[A(n,s)=\frac{2^{-n-1}}{1-2^{1-s}} \sum_{k=0}^n\frac{(-n)_k}{k!}\,\frac{1}{(k+1)^s},\]
where $(a)_k=\g(a+k)/\g(a)$ is the Pochhammer symbol. Using the standard result
\[\frac{1}{(k+1)^s}=\frac{1}{\g(s)} \int_0^\infty e^{-(k+1)w} w^{s-1} dw\qquad (\sigma>0),\]
we find upon inversion of the order of summation and integration that
\begin{eqnarray}
A(n,s)&=&\frac{2^{-n-1}}{(1-2^{1-s}) \g(s)}\int_0^\infty e^{-w} w^{s-1} \sum_{k=0}^n \frac{(-n)_k e^{-kw}}{k!} dw\nonumber\\
&=&\frac{2^{-n-1}}{(1-2^{1-s}) \g(s)}\int_0^\infty e^{-w}(1-e^{-w})^n w^{s-1}dw.\label{e21}
\end{eqnarray}
This result, which holds by analytic continuation for $\sigma>-n$, was established in \cite[Eq.~(17)]{Jerby} by a more elaborate derivation.

With $s=\sigma+it$, $t>0$, we now set $t=an$ as $n\to\infty$, where $a>0$ is a finite constant. We note that the factor controlling convergence of the integral at infinity is $e^{-w}$.  Thus, we put 
\bee\label{e22}
I(n,a;\sigma):=\int_0^\infty e^{-w}(1-e^{-w})^n w^{s-1}dw=\int_0^\infty w^{\sigma-1} e^{n\psi(w,n)}dw,
\ee
where
\[\psi(w,n):=\log (1-e^{-w})+ia \log\,w-\frac{w}{n}.\]
It is clear that the saddle points of the phase function $\psi(w,n)$, given by those points where $\psi'(w,n) 
=0$, will depend on the large parameter $n$; when $|w|$ is not large this dependence will be weak, although when $|w|$ is larger this dependence will become more significant.
It is necessary to include the factor $e^{-w}$ in $\psi(w,n)$ in order for the paths of steepest descent to connect with the point at infinity in $\Re (w)>0$; see \cite[p.~42]{ETC}.

The saddle points are given by the roots of the equation
\bee\label{e23}
\frac{1}{e^w-1}+\frac{ia}{w}-\frac{1}{n}=0
\ee
and for $t>0$ are found to be situated\footnote{We do not need to consider the saddles in $\Im (w)<0$.} in an infinite string approximately parallel to the positive imaginary $w$-axis.
We label the saddles $w_k=x_k+iy_k$, $k=1, 2, \ldots$ and observe from (\ref{e23}) that
\bee\label{e24}
x_k\simeq \log\,((2k-1)\pi/a),\qquad y_k\simeq (2k-1)\pi
\ee
for the low-lying saddles. The values of $w_k$ were computed from (\ref{e23}) using the FindRoot command in {\it Mathematica} with (\ref{e24}) as starting values. The steepest descent paths through these saddles terminate at the points $2\pi ki$, $k=0, 1, 2, \ldots\ $, which are logarithmic singularities of the phase function. Branch cuts emanating from these singularities can be taken along the horizontal lines $(-\infty+2\pi ki,2\pi ki]$. A typical example of these steepest descent paths is shown in Fig.~1 for the case $n=20$, $a=1$, where the contributory saddles $w_k$ correspond to $1\leq k\leq4$. The steepest descent path through $w_4$ has one endpoint at $w=6\pi i$ with the other passing to infinity in $\Re (w)>0$. The steepest ascent paths through $w_k$ ($1\leq k\leq m-1$) all encircle the origin and pass to infinity in $\Re (w)<0$.
\begin{figure}[th]
	\begin{center}\includegraphics[width=0.6\textwidth]{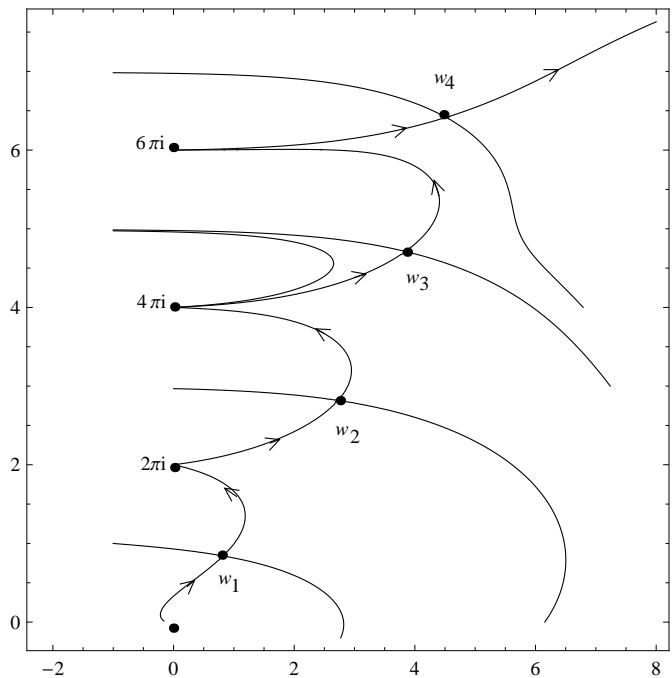}
	
\caption{\small{The paths of steepest descent through the contributory saddles $w_k$, $1\leq k\leq4$ when $n=20$, $a=1$. The arrows indicate the direction of integration. The vertical scale represents multiples of $\pi i$.}}
\end{center}
\end{figure}
It is found by observation that the steepest descent path that passes through the saddle labelled by $k=m$ and thence
goes to infinity in $\Re (w)>0$ corresponds to when $a\simeq y_m/n$; that is, from (\ref{e24}),
\bee\label{e25}
m\simeq \bl\lceil\frac{t}{2\pi}+\frac{1}{2}\br\rceil.
\ee

The method of steepest descents applied to the integral in (\ref{e22}) shows that the contribution to $I(n,a;\sigma)$ from the $k$th saddle is given by (see, for example, \cite[p.~47]{DLMF}, \cite[p.~13]{PHad})
\bee\label{e26}
I_k(n,a;\sigma)\sim i\,\sqrt{\frac{2\pi}n }\,\frac{ e^{-w_k}(1-e^{-w_k})^n}{\sqrt{\psi''(w_k,n)}}\,w_k^{s-1} \sum_{j=0}^\infty \frac{c_j}{n^j}\,\frac{\g(j+\fs)}{\g(\fs)},
\ee
where
\[\psi''(w_k,n)=\frac{1}{e^{w_k}-1}-\frac{ia}{w_k}=\frac{1}{1-e^{-w_k}}\bl(\frac{ia}{w_k}-\frac{1}{n}\br)-\frac{ia}{w_k^2}.\]
The first three coefficients $c_j$ have the form \cite[pp.~13--14]{PHad}
\[c_0=1,\quad c_1=\frac{-1}{2\psi''(w_k,n)}\bl\{2F_2-2\Psi_3 F_1+\frac{5}{6}\Psi_3^2-\frac{1}{2}\Psi_4\br\},\]
\[c_2=\frac{1}{(2\psi''(w_k,n))^2}\bl\{\frac{2}{3}F_4-\frac{20}{9}\Psi_3F_3+\frac{5}{3}\bl(\frac{7}{3}\Psi_3^2-\Psi_4\br)F_2-\frac{35}{9}\bl(\Psi_3^3-\Psi_3\Psi_4+\frac{6}{35}\Psi_5\br)F_1\]
\[+\frac{35}{9}\bl(\frac{11}{24}\Psi_3^4-\frac{3}{4}\bl(\Psi_3^2-\frac{1}{6}\Psi_4\br)\Psi_4+\frac{1}{5}\Psi_3\Psi_5-\frac{1}{35}\Psi_6\br)\br\},\]
where
\[\Psi_j:=\frac{\psi^{(j)}(w;n)}{\psi''(w,n)},\qquad F_j:=\frac{f^{(j)}(w)}{f(w)},\]
with $f(w):=w^{\sigma-1}$ and the derivatives evaluated at the saddle point $w=w_k$.

Application of Cauchy's theorem then shows that the integration path $[0,\infty)$ can be deformed to pass over the steepest descent paths through the saddles $1\leq k\leq m$ to yield
\bee\label{e27}
A(n,s)=\frac{2^{-n-1}}{(1-2^{1-s}) \g(s)}\,\sum_{k=0}^m I_k(n,a;\sigma)
\ee
as $n\to\infty$, where $I_k(n,a;\sa)$ has the asymptotic expansion given by (\ref{e26}).

\vspace{0.6cm}

\begin{center}
{\bf 3.\ Numerical verification of the expansion}
\end{center}
\setcounter{section}{3}
\setcounter{equation}{0}
\renewcommand{\theequation}{\arabic{section}.\arabic{equation}}
In this section we present some numerical examples to illustrate the validity of the result in (\ref{e27}). An example of the contributory saddles $w_k$ for $n=40$, $a=1$ and $n=20$, $a=2$ (both cases corresponding to $m=7$) are shown in Table 1.
In Tables 2 and 3, we present values of $A(n,s)$ computed from (\ref{e12}) compared with the asymptotic result from (\ref{e27}) with $j\leq 2$ in the particular case $\sigma=\fs$ (on the critical line $\Re (s)=\fs$). In each case, the value of $m$ corresponding to the last contributory saddle with steepest descent path passing from the singularity $w=2\pi(m-1)i$ to infinity in $\Re (w)>0$ is indicated.

\begin{table}[h]
\caption{\footnotesize{Values of the contributory saddles $w_k$ (to 6dp) for two cases of $n$ and $a$ (both corresponding to $m=7$).}}
\begin{center}
\begin{tabular}{|c|l|l|}
\hline
&&\\[-0.3cm]
\mcol{1}{|c|}{$k$} & \mcol{1}{c|}{$n=40,\ a=1$} & \mcol{1}{c|}{$n=20,\ a=2$}\\
\hline
&&\\[-0.3cm]
1 & $-0.213894 + 3.299584i$  & $+0.735036 + 2.723878i$ \\
2 & $+1.619139 + 9.152549i$  & $+2.410605 + 9.057147i$ \\
3 & $+2.458648 + 15.428395i$ & $+3.191141 + 15.360866i$ \\
4 & $+3.117553 + 21.666246i$ & $+3.823744 + 21.602927i$ \\
5 & $+3.765495 + 27.830032i$ & $+4.448382 + 27.761143i$ \\
6 & $+4.448382 + 27.761143i$ & $+5.121844 + 33.718312i$ \\
7 & $+5.045810 + 39.277859i$ & $+5.632058 + 39.255291i$ \\
\hline
\end{tabular}
\end{center}
\end{table}
\begin{table}[h]
\caption{\footnotesize{Values of $A(n,s)$ and the asymptotic value from (\ref{e27}) when $n=20$ and $\sigma=\fs$.}}
\begin{center}
\begin{tabular}{|cc|l|l|}
\hline
&&&\\[-0.3cm]
\mcol{1}{|c}{$a$} & \mcol{1}{c|}{$m$} & \mcol{1}{c|}{$A(n,s)$} & \mcol{1}{c|}{\mbox{Asymptotic value}}\\
\hline
&&&\\[-0.3cm]
0.50 & 2 & $-0.0002394854+0.0000979486i$ & $-0.00023983+0.00009811i$\\
0.75 & 3 & $-0.0013656997-0.0009979383i$ & $-0.00136554-0.00099839i$\\
0.80 & 3 & $-0.0026415717+0.0020871724i$ & $-0.00264151+0.00208667i$\\
1.00 & 4 & $+0.0086008223-0.0117220182i$ & $+0.00860160-0.01720826i$\\
1.50 & 5 & $-0.0511931929+0.0054038870i$ & $-0.05119219+0.00540340i$\\
2.00 & 7 & $-0.0085839350-0.0372653861i$ & $-0.00858386-0.03726493i$\\
5.00 & 17& $-0.1462531266-0.0449764455i$ & $-0.14625160-0.04497750i$\\
\hline
\end{tabular}
\end{center}
\end{table}
\begin{table}[h]
\caption{\footnotesize{Values of $A(n,s)$ and the asymptotic value from (\ref{e27}) when $n=50$ and $\sigma=\fs$.}}
\begin{center}
\begin{tabular}{|cc|l|l|}
\hline
&&&\\[-0.3cm]
\mcol{1}{|c}{$a$} & \mcol{1}{c|}{$m$} & \mcol{1}{c|}{$A(n,s)$} & \mcol{1}{c|}{\mbox{Asymptotic value}}\\
\hline
&&&\\[-0.3cm]
0.80 & 7 & $+0.0000234378+0.0000433293i$ & $+0.0000234374+0.0000433292i$\\
1.00 & 9 & $+0.0004150615-0.0009392525i$ & $+0.0004150622-0.0009392487i$\\
1.50 & 13& $-0.0353214881-0.0050091223i$ & $-0.0353214525-0.0050091204i$\\
2.00 & 17& $+0.0460334465+0.0392889898i$ & $+0.0460334317+0.0392889689i$\\
4.00 & 33& $+0.0242455885-0.0183724506i$ & $+0.0242455076-0.0183724384i$\\
5.00 & 41& $+0.0188678860+0.0014542050i$ & $+0.0188678811+0.0014542105i$\\
\hline
\end{tabular}
\end{center}
\end{table}

It will be observed that as $n$ and $t$ become large the number of contributory saddles increases. In Fig.~2(a) we show the variation of the absolute value of the contribution from the $k$th saddle
\[{\hat I}_k:=\bl|\frac{2^{-n-1}}{(1-2^{1-s}) \g(s)}\,I_k(n,a;\sa)\br|\]
by plotting $\log_{10}{\hat I}_k$ against $k$ for $n=50$ and $a=5$, which reveals a reasonably strong decay in the saddle contributions. As a consequence, it is feasible to truncate the number of contributory saddles considered depending on the accuracy required.
\begin{figure}[th]
	\begin{center}	{\tiny($a$)}\ \includegraphics[width=0.4\textwidth]{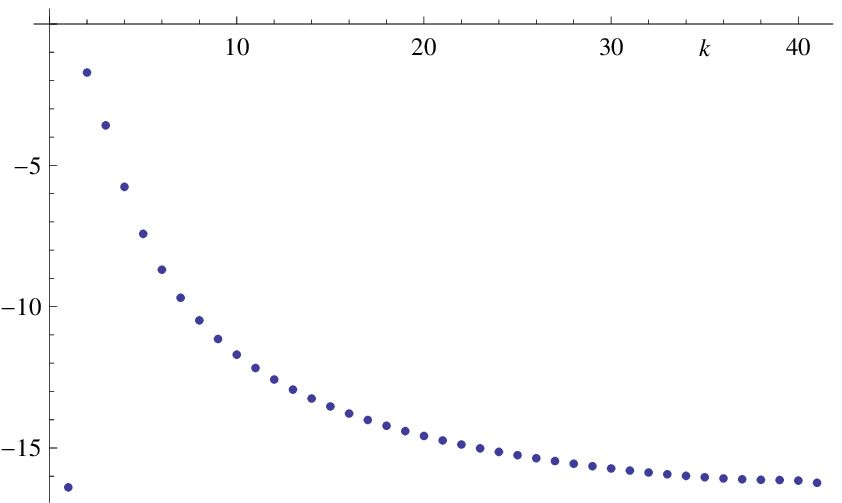}
	\qquad
	{\tiny($b$)}\ \includegraphics[width=0.38\textwidth]{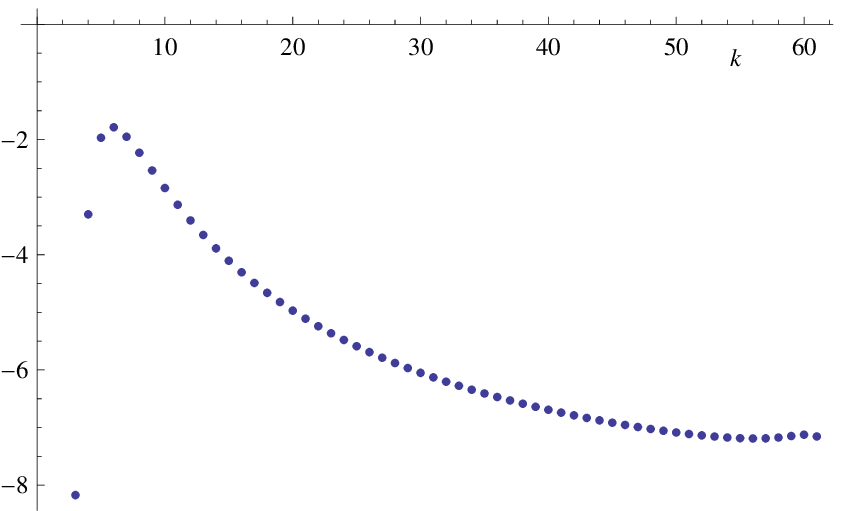} 
	
\caption{\small{Plots of $\log_{10}{\hat I}_k$ against $k$ showing the variation in the absolute value of the saddle contributions for $\sa=\fs$ when (a) $n=50$, $a=5$ ($1\leq k\leq m$) and (b) $n=20$, $a=6\pi$ ($3\leq k\leq m$).}}
\end{center}
\end{figure}

It is worth noting that it is possible for the saddles $w_{m-1}$ and $w_m$ to connect resulting in a Stokes phenomenon. This arises when the steepest descent path through $w_{m-1}$, instead of terminating at the point $w=2\pi (m-1)i$, connects with the saddle $w_m$. This occurs for example, when $n=5$ and $a\doteq 6.032$ ($t\doteq 30.160$).  Since this can only arise for the last two saddles in the contributory sequence, where the contribution to $A(n,s)$ is small, we do not consider this complication here.

The values of the parameter $a$ considered by Jerby in \cite{Jerby} are given by $a=\pi N$, where integer $N\geq 1$. For these larger values of $a$ it is found that the first $N$ saddles lie in $\Re (w)<0$ (when $n$ is large), as observed in \cite{Jerby}. It is found by observation that the saddles corresponding to $1\leq k\leq k^*-1$, $k^*=\lfloor (N+1)/2\rfloor$ are non-contributory as these saddles cannot connect with those for which the steepest descent path ultimately passes to infinity in $\Re (w)>0$. The steepest descent path from the origin (which spirals out of the origin over different Riemann sheets) that connects with the first saddle $w_{k^*}$ then proceeds over the remaining contributory saddles with $k^*+1\leq k\leq m$ in a similar manner to that shown in Fig.~1. An example of the behaviour of the absolute value of the saddle contributions when $n=20$ and $a=6\pi$ is shown in Fig.~2(b). 
In Table 4 we show values of $A(n,s)$ when $a=\pi N$ compared with the asymptotic values obtained from (\ref{e27}) with $j\leq 2$ in the particular case $\sigma=\fs$.
%
\begin{table}[h]
\caption{\footnotesize{Values of $A(n,s)$ and the asymptotic value from (\ref{e27}) when $n=30$, $a=\pi N$ and $\sigma=\fs$.}}
\begin{center}
\begin{tabular}{|cc|l|l|}
\hline
&&&\\[-0.3cm]
\mcol{1}{|c}{$N$} & \mcol{1}{c|}{$m$} & \mcol{1}{c|}{$A(n,s)$} & \mcol{1}{c|}{\mbox{Asymptotic value}}\\
\hline
&&&\\[-0.3cm]
1 & 16 & $+0.0021433151+0.0011784556i$ & $+0.0021433011+0.0011784496i$\\
2 & 31 & $+0.0120051627+0.0069585493i$ & $+0.0120052138+0.0069585241i$\\
3 & 46 & $-0.0288262956+0.0163914511i$ & $-0.0288262658+0.0163913977i$\\
4 & 61 & $+0.0053628619+0.0257175197i$ & $+0.0053628513+0.0257174689i$\\
5 & 76 & $+0.0929962033+0.0664340984i$ & $+0.0929959750+0.0664339537i$\\
\hline
\end{tabular}
\end{center}
\end{table}

\begin{figure}[h]
	\begin{center}	{\tiny($a$)}\ \includegraphics[width=0.5\textwidth]{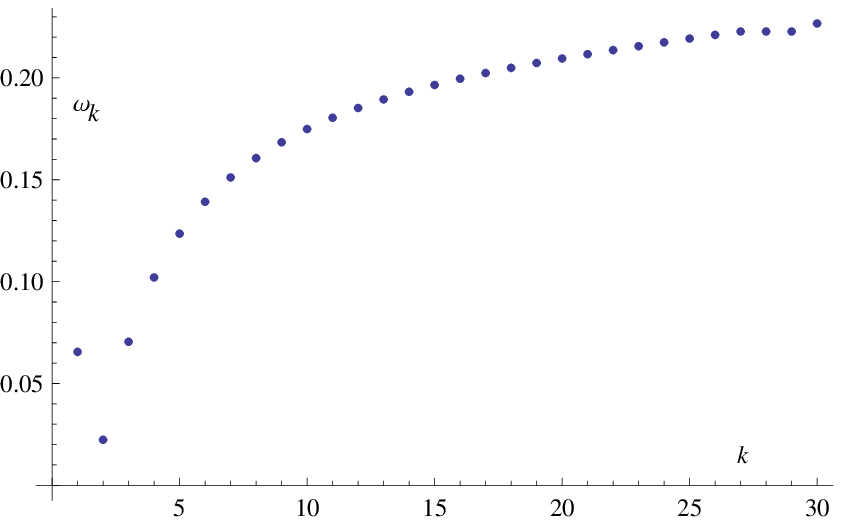}
	
\caption{\small{Plot of $\om_k$ against $k$ when $n=50$, $a=\pi$ ($m=31$) and $\sa=\fs$.}}
\end{center}
\end{figure}

By inspection, it is found that when $a=\pi N$ (and $n$ is sufficiently large) the maximum saddle point contribution\footnote{The maximum saddle point contribution corresponds to $k=1$ when $0<a\,\ltwid \, 0.88\pi$.} (in absolute value) corresponds to $k=N$ when $N>1$. Referring to Fig.~2, the maximum saddle contribution occurs for $k=2$ in the first figure and $k=6$ in the second figure. In the case of Fig.~2(a), the saddle $w_2$ produces the dominant contribution to $A(n,s)$ and we find ${\hat I}_2\doteq 0.019205$ with $|A(n,s)|\doteq 0.018924$, thus confirming that the main contribution arises from the saddle corresponding to $k=2$. In Fig.~2(b) there are neighbouring saddles yielding a comparable magnitude to that corresponding to $w_6$, which will result in a significant cancellation between the real and imaginary parts; in this case we find ${\hat I}_6\doteq 0.01634$ and $|A(n,s)|\doteq 0.03653$.

A primary concern in \cite{Jerby} was the demonstration that $A(n,s)$ is exponentially small as $t\to+\infty$.
From (\ref{e26}) and (\ref{e27}), the order of ${\hat I}_k$ is controlled by
$O(e^{-na\om_k})$, or equivalently $O(e^{-\om_k t})$, where
\[\om_k:=\theta_k-\frac{1}{2}\pi-\frac{1}{a} \log\bl|\frac{1-e^{-w_k}}{2}\br|+\frac{\Re (w_k)}{t},\qquad \theta_k:=\arg\,w_k,\]
and we have used the fact that $|\g(\sa+it)|\sim \sqrt{2\pi}\,t^{\sa-1/2} e^{-\pi t/2}$ as $t\to+\infty$ \cite[(5.11.9)]{DLMF}. It is found that $\om_k>0$ for all contributory saddles; a typical plot of $\om_k$ against $k$
in the case $n=50$, $a=\pi$ ($m=31$) is shown in Fig.~3. The value of $\om_k$ for the saddle with $k=2$
is approximately\footnote{In the absence of the term involving $t^{-1}$ in $\om_k$ the minimum value is $\om_2\doteq0.01773$, as found in \cite{Jerby}.} 0.02235, with the values for $k\geq 3$ increasing monotonically. The small value of $\om_2$ would make it difficult to estimate the order of $A(n,s)$ by simply taking the smallest value ${\hat I}_2$ and multiplying by the number of contributory saddles of $O(t/2\pi)$.


\vspace{0.6cm}

\begin{thebibliography}{99}
\footnotesize{

\bibitem{ETC}
E.T. Copson, {\it Asymptotic Expansions}, Cambridge University Press, Cambridge, 1965.

\bibitem{Hasse}
H. Hasse, Ein Summierungsverfahren fur die Riemannsche $\zeta$-Reihe, Math. Z. {\bf 32} (1930) 458--464.

\bibitem{Ivic}
A. Ivi\'c, {\it The Theory of the Riemann-zeta Function with Applications}, J. Wiley \& Sons, New York, 1985.

\bibitem{Jerby}
Y. Jerby, An approximate functional equation for the Riemann zeta function with exponentially decaying error, J. Approx. Theory {\bf 265} (2021) 105551.

\bibitem{DLMF}
F.W.J. Olver, D.W. Lozier, R.F. Boisvert and C.W. Clark (eds.),    
{\it NIST Handbook of Mathematical Functions}, Cambridge University Press, Cambridge, 2010.

\bibitem{PHad}
R.B. Paris, {\it Hadamard Expansions and Hyperasymptotic Evaluation}, Encyclopedia of Mathematics and its Applications Vol. 141, Cambridge University Press, Cambridge, 2011.

\bibitem{Sond}
J. Sondow, Analytic continuation of Riemann's zeta function and values at negative integers via Euler's transformation of series, Proc. Amer. Math. Soc. {\bf 120} (1994) 421--425.



}
\end{thebibliography}
\end{document}